\documentclass{amsart}

\usepackage{amsmath, amssymb, latexsym, amscd, stmaryrd}
\usepackage{tikz,tikz-cd}

\makeatletter
\DeclareFontFamily{OMX}{MnSymbolE}{}
\DeclareSymbolFont{MnLargeSymbols}{OMX}{MnSymbolE}{m}{n}
\SetSymbolFont{MnLargeSymbols}{bold}{OMX}{MnSymbolE}{b}{n}
\DeclareFontShape{OMX}{MnSymbolE}{m}{n}{
    <-6>  MnSymbolE5
   <6-7>  MnSymbolE6
   <7-8>  MnSymbolE7
   <8-9>  MnSymbolE8
   <9-10> MnSymbolE9
  <10-12> MnSymbolE10
  <12->   MnSymbolE12
}{}
\DeclareFontShape{OMX}{MnSymbolE}{b}{n}{
    <-6>  MnSymbolE-Bold5
   <6-7>  MnSymbolE-Bold6
   <7-8>  MnSymbolE-Bold7
   <8-9>  MnSymbolE-Bold8
   <9-10> MnSymbolE-Bold9
  <10-12> MnSymbolE-Bold10
  <12->   MnSymbolE-Bold12
}{}

\let\llangle\@undefined
\let\rrangle\@undefined
\DeclareMathDelimiter{\llangle}{\mathopen}%
                     {MnLargeSymbols}{'164}{MnLargeSymbols}{'164}
\DeclareMathDelimiter{\rrangle}{\mathclose}%
                     {MnLargeSymbols}{'171}{MnLargeSymbols}{'171}
\makeatother

\newcommand{\bC}{\mathbb{C}}

\newcommand{\bL}{\mathbb{L}}
\newcommand{\bP}{\mathbb{P}}

\newcommand{\bR}{\mathbb{R}}
\newcommand{\bT}{\mathbb{T}}
\newcommand{\bZ}{\mathbb{Z}}

\newcommand{\cB}{\mathcal{B}}

\newcommand{\cF}{\mathcal{F}}
\newcommand{\cL}{\mathcal{L}}
\newcommand{\cO}{\mathcal{O}}
\newcommand{\cM}{\mathcal{M}}
\newcommand{\cX}{\mathcal{X}}

\newcommand{\IX}{\mathcal{IX}}

\newcommand{\hX}{\hat{X}}

\newcommand{\tL}{\widetilde{L}}

\newcommand{\tmu}{\widetilde{\mu}}
\newcommand{\tW}{\widetilde{W}}

\newcommand{\Mbar}{\overline{\cM} }

\newcommand{\Aut}{\mathrm{Aut}}
\newcommand{\vir}{\mathrm{vir}}
\newcommand{\ev}{\mathrm{ev}}

\newcommand{\age}{\mathrm{age}}
\newcommand{\CR}{\mathrm{CR}}
\newcommand{\Hom}{\mathrm{Hom}}

\newcommand{\Crit}{ {\mathrm{Crit}} }
\newcommand{\Res}{ \mathrm{Res} }
\newcommand{\Spec}{\mathrm{Spec}}
\newcommand{\BSi}{ \mathrm{Box}(\Si)}
\newcommand{\Int}{\mathrm{Int}}
\newcommand{\Area}{\mathrm{Area}}
\newcommand{\SYZ}{\mathrm{SYZ}}

\newcommand{\fl}{\mathfrak{l}}

\newcommand{\Si}{ {\Sigma} }
\newcommand{\si}{\sigma}

\newcommand{\vmu}{{\vec{\mu}}}

\newcommand{\ww}{\mathsf{w}}
\newcommand{\one}{\mathbf{1}}
\newcommand{\btau}{\boldsymbol{\tau}}

\newcommand{\vGa}{\vec{\Gamma}}

\newcommand{\fg}{\mathfrak{g}}
\newcommand{\fn}{\mathfrak{n}}

\newcommand{\Cbar}{\overline{C}}

\newtheorem{theorem}{Theorem}

\newtheorem{conjecture}[theorem]{Conjecture}
\theoremstyle{definition}

\theoremstyle{remark}

\numberwithin{equation}{section}

\begin{document}

\title[SYZ mirror symmetry \& BKMP remodeling conjecture]{The SYZ mirror symmetry and the BKMP remodeling conjecture}

\author{Bohan Fang}

\address{Bohan Fang\\
Beijing International Center for Mathematical Research\\
 Peking University, 5 Yiheyuan Road, Beijing 100871, China} 

\email{bohanfang@gmail.com}

\author{Chiu-Chu Melissa Liu}

\address{Chiu-Chu Melissa Liu\\
Department of Mathematics, Columbia University\\
 2990 Broadway, New York, NY 10027, USA}

\email{ccliu@math.columbia.edu}

\author{Zhengyu Zong}

\address{Zhengyu Zong\\
Yau Mathematical Sciences Center, Tsinghua University\\
Jin Chun Yuan West Build., Tsinghua University,  Beijing 100084, China}

\email{zyzong@math.tsinghua.edu.cn}

\begin{abstract}
The Remodeling Conjecture proposed by Bouchard-Klemm-Mari\~{n}o-Pasquetti (BKMP) relates the A-model open and closed
topological string amplitudes (open and closed Gromov-Witten invariants) of a symplectic toric Calabi-Yau 3-fold to
Eynard-Orantin invariants of its mirror curve. The Remodeling Conjecture can be viewed as a version of
all genus open-closed mirror symmetry. The SYZ conjecture explains mirror symmetry as $T$-duality. After a brief review on SYZ mirror symmetry
and mirrors of symplectic toric Calabi-Yau 
3-orbifolds, we give a non-technical exposition of our results on the Remodeling Conjecture for symplectic toric Calabi-Yau 3-orbifolds. In the end, we apply SYZ mirror symmetry to obtain the descendent version of the all genus mirror symmetry for toric Calabi-Yau 3-orbifolds.
\end{abstract}

\maketitle

\section{Introduction}

\subsection{The SYZ conjecture}
Mirror symmetry relates the A-model on a Calabi-Yau $n$-fold $\cX$, defined
by the symplectic structure on $\cX$, to the B-model on
a mirror Calabi-Yau $n$-fold $\check{X}$, defined by the complex
structure on $\check{X}$. Strominger-Yau-Zaslow proposed that
Mirror Symmetry is T-duality \cite{SYZ} in the following sense.
There exist $\pi: \cX\to B$ and $\check{\pi}:\check{\cX}\to B$, where
$\dim_\bR B=n$,  such
that over a dense open subset $U\subset B$,
$Y:=\pi^{-1}(U)\to U$ and $\check{Y}:=\check{\pi}^{-1}(U)\to U$ are dual special Lagrangian
$n$-torus fibrations.
In the semi-flat case, $U=B$,  $Y\cong T^*B/\Lambda$
as a symplectic manifold and $\check{Y}\cong TB/\Lambda^\vee$
as a complex manifold, where $\Lambda$ and $\Lambda^\vee$ are dual lattices,
and for each $b\in B$, $\pi^{-1}(b)=T^*_b B/\Lambda$ and
$\check{\pi}^{-1}(b) = T_b B/\Lambda^\vee$. In the general case, there is quantum
correction to the complex structure on $\check{Y}$ coming from holomorphic disks
in $\cX$ bounded by  Lagrangian $n$-tori $\pi^{-1}(b)$, $b\in U$.

Mirror symmetry has been extended to certain non-Calabi-Yau manifolds.
When $\cX$ is Fano (or more generally semi-Fano), the mirror is a Landau-Ginzburg model, which can also be constructed by SYZ transformation.

\subsection{SYZ mirror symmetry for toric Calabi-Yau 3-folds}
Let $(\cX,\omega)$ be a symplectic toric Calabi-Yau 3-manifold/orbifold,
where $\omega$ is the symplectic form. There are two families of mirrors,
both of which can be constructed by the SYZ transformation.

\subsubsection*{Landau-Ginzburg mirror}
The mirror B-model to the toric Calabi-Yau 3-orbifold $(\cX,\omega)$ is a
$3$-dimensional Landau-Ginzburg model on $(\bC^*)^3$ given by the superpotential
$$
W=H(X,Y,q)Z.
$$
The Calabi-Yau condition ensures that $W$ is in this form. The
complex parameter $q$ is related to the K\"ahler parameter of $\cX$ by
the mirror map.
The Landau-Ginzburg mirror $( (\bC^*)^3,W)$ can be constructed
by applying SYZ transformation to $\mu: \cX\to \Delta\subset \bR^3$,
where $\mu$ is the moment map of the Hamiltonian $U(1)^3$-action
on $(\cX,\omega)$, and $\Delta=\mu(\cX)$ is the moment polyhedron.

\subsubsection*{Hori-Vafa mirror}
By Hori-Vafa \cite{HV}, the mirror of $(\cX,\Omega)$ is a non-compact Calabi-Yau 3-fold
$(\check{\cX},\Omega)$, where
$$
\check{\cX}=\{ (u,v,X,Y): u,v\in \bC, X,Y\in \bC^*, uv=H(X,Y)\}
$$
is a hypersurface in $\bC^2\times (\bC^*)^2$, and
$$
\Omega =\Res_{ \check{\cX} }\left(\frac{1}{uv-H(X,Y)}du\wedge dv\wedge\frac{dX}{X}\wedge \frac{dY}{Y}\right)
$$
is a holomorphic 3-form on $\check{\cX}$. The Hori-Vafa mirror
$(\check{\cX},\Omega)$ can be constructed by applying SYZ transformation
to the Gross fibration \cite{CLL, CCLT}.

\bigskip

The two ``equivalent'' mirrors come with no surprise since as a toric
variety $\cX$ should have a Landau-Ginzburg mirror, while as a Calabi-Yau $3$-fold
Hori-Vafa  showed that they both could be reduced to a mirror curve
$$
C_q=\{ (X,Y)\in (\bC^*)^2: H(X,Y,q)=0\} \subset (\bC^*)^2.
$$
In particular
\begin{equation}\label{eqn:OSintegral}
\int_{\Gamma} e^{-W}\frac{dXdYdZ}{XYZ} = \int_{\widetilde \Gamma}
\Omega = \int_{\gamma} ydx,
\end{equation}
where $x=-\log X$, $y=-\log Y$,  and Lagrangian cycles $\Gamma\subset (\bC^*)^3$,
$\widetilde\Gamma\subset \check{\cX}$, $\gamma\subset C_q$ are related by a series of
dimensional reductions.

\subsection{The BKMP remodeling conjecture}
It is usually difficult to obtain higher genus
invariants of the B-model. A standard way is to
apply the BCOV holomorphic anomaly equations \cite{BCOV93}. These equations do not
have unique solutions (holomorphic ambiguity). One needs to fix the boundary
conditions via the input from mirror symmetry, like Gromov-Witten
invariants in low degrees in the large radius limit point, and to utilize the so-called
``gap conditions'' at the conifold point. To mathematically prove
the all genus mirror symmetry under this B-model approach is currently
beyond reach since the A-side theory of BCOV (and the A-side theory at
the conifold point) is still a mystery to mathematicians.

The Eynard-Orantin topological recursion is an algorithm which produces higher genus
invariants for a spectral curve \cite{EO07}. Applying the Eynard-Orantin
topological recursion to the mirror curve $C_q$ of a symplectic toric Calabi-Yau 3-orbifold
$(\cX,\omega)$, we obtain a version of the B-model, related but not a priori the same as the BCOV theory on the Hori-Vafa mirror $(\check{\cX},\Omega)$ of $(\cX,\omega)$.
The Bouchard-Klemm-Mari\~{n}o-Pasquetti (BKMP) remodeling
conjecture \cite{BKMP09, BKMP10} says that all genus B-model topological strings on $(\check{\cX},\Omega)$ are essentially Eynard-Orantin invariants \cite{EO07} of
$C_q$. Using mirror symmetry, BKMP relates the Eynard-Orantin invariant $\omega_{g,n}$
of the mirror curve to a generating function $F^{\cX,\cL}_{g,n}$ of open
GW invariants (A-model topological open string amplitudes) counting holomorphic
maps from bordered Riemann surfaces with $g$ handles and $n$ holes to $\cX$ with
boundaries in an Aganagic-Vafa Lagrangian brane $\cL$. The
correspondence between $\omega_{g,n}$ and $F^{\cX,\cL}_{g,n}$ can be
extend to the case when $n=0$. In this case, $\omega_{g,0}$ is defined
to be the free energy and the A-model potential becomes the closed
Gromov-Witten potential $F^\cX_g$. Therefore, the remodeling
conjecture gives us an all genus open-closed mirror symmetry for toric
Calabi-Yau 3-orbifolds. The three equivalent B-models mirror to $\cX$ are
illustrated as below.
$$
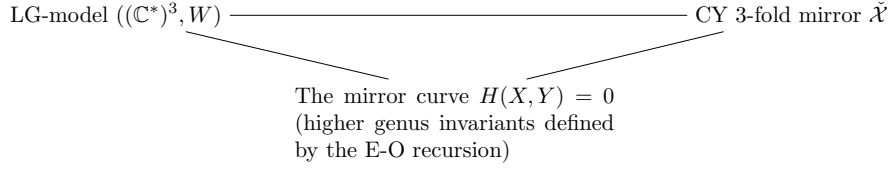
\begin{tikzpicture}[baseline= (a).base]
\node[scale=.85] (a) at (0,0){
\begin{tikzcd}
\text{LG-model $((\bC^*)^3,W)$} \arrow[-]{rr}\arrow[-]{rd}
&
& \text{CY $3$-fold mirror $\check \cX$} \arrow[-]{ld}\\
& \text{\parbox{5cm}{\bigskip \bigskip \bigskip The mirror curve $H(X,Y)=0$ (higher genus
    invariants defined by the E-O recursion)} }
\end{tikzcd}};
\end{tikzpicture}
$$
\bigskip

The BKMP remodeling conjecture was proved for $\bC^3$ at all genus $g$ independently by L. Chen \cite{Ch09}
and J. Zhou \cite{Zh09} in the $n>0$ case (open string sector), and
by Bouchard-Catuneanu-Marchal-Su{\l}kowski \cite{BCMS} in the $n=0$ case (closed string sector).
In 2012, Eynard-Orantin provided a proof of the BKMP remodeling conjecture for all symplectic smooth toric Calabi-Yau 3-folds \cite{EO12}.

The SYZ T-duality transformation \cite{SYZ, AP, LYZ} associates a
coherent sheaf $\cF$ on $\cX$ to a Lagrangian cycle $\SYZ(\cF)\subset (\bC^*)^3$.
A coherent sheaf on $\cX$ is admissible if there exists $\gamma(\cF) \subset C_q$ such that
$$
\int_{\SYZ(\cF)}e^{\frac{W +x+ fy}{z}}\frac{dX dY dZ}{XYZ} =\int_{\gamma(\cF)}e^{\frac{x+fy}{z}} ydx
$$
and $x+fy$ is bounded below on $\gamma(\cF)$. Here $f$ is
an integer and $z$ is a negative real number, so that the
integral on the right hand side converges.
The remodeling conjecture has a descendant
version: given $n$ admissible coherent sheaves $\cF_1,\ldots,\cF_n$ on $\cX$,
the Laplace transform of $\omega_{g,n}$
along $\gamma(\cF_1),\ldots,\gamma(\cF_n)$ is a generating function of genus $g$ descendant
Gromov-Witten invariants with $n$ insertions
$\kappa(\cF_1),\ldots, \kappa(\cF_n)$, where $\kappa(\cF_i)$ is
the so-called Gamma class of $\cF_i$.

In the rest of this paper, we will give a non-technical exposition of our
results on the remodeling conjecture for toric Calabi-Yau 3-orbifolds \cite{FLZ13, FLZ15}, which is a version of all genus open-closed mirror symmetry. We will also
discuss the all genus mirror symmetry of free energies and descendant potentials.

\section{Toric Calabi-Yau 3-orbifolds and their mirror curves}

\subsection{Toric Calabi-Yau 3-orbifolds}
A Calabi-Yau 3-fold $X$ is toric if it contains the algebraic torus
$\bT=(\bC^*)^3$ as a Zariski dense open subset, and the action of
$\bT$ on itself extends to $X$. All Calabi-Yau 3-folds are non-compact.
There is a rank 2 subtorus $\bT' \subset \bT$
which acts trivially on the canonical line bundle of $X$.
We call $\bT'$ the Calabi-Yau torus. Then $\bT\cong \bT'\times \bC^*$. Let $\bT'_\bR\cong U(1)^2$ be the
maximal compact subgroup of $\bT'$.

Let $M'=\Hom(\bT',\bC^*)\cong \bZ^2$  and $N'=\Hom(\bC^*,\bT')$ be the character lattice
and the cocharacter lattice of $\bT'$, respectively. Then $M'$ and $N'$
are dual lattices.  Let $X_\Sigma$ be a toric Calabi-Yau 3-fold defined by a
simplicial fan $\Si \subset  N'_\bR\times \bR$,
where $N'_\bR:=N'\otimes_{\bZ}\bR\cong \bR^2$ can be identified with
the Lie algebra of $\bT_\bR'$.
Then $X_\Si$ has at most quotient singularities. We assume that
$X_\Si$ is semi-projective, i.e., $X_\Si$ contains at least one $\bT$
fixed point, and $X_\Si$ is projective over its affinization
$X_0:=\Spec H^0(X_\Si,\cO_{X_\Si})$.  Then the support of the fan $\Si$ is
a strongly convex rational polyhedral cone $\si_0\subset N'_\bR \times \bR \cong \bR^3$,
and $X_0$ is the affine toric variety defined by the 3-dimensional cone $\si_0$.
There exists a convex polytope $P \subset N'_\bR\cong \bR^2$ with vertices in the lattice $N'\cong \bZ^2$, such
that $\si_0$ is the cone over $P\times \{1\}\subset N'_\bR\times \bR $, i.e.
$\si_0=\{ (tx, ty, t): (x,y)\in P, t\in [0,\infty)\}$.
The fan $\Si$ determines a triangulation of $P$:
the 1-dimensional, 2-dimensional, and 3-dimensional cones
in $\Si$ are in one-to-one correspondence with
the vertices, edges, and faces of the triangulation of $P$, respectively.
This triangulation of $P$ is known as the toric diagram or the dual graph of the simplicial
toric Calabi-Yau 3-fold $X_\Si$.

Let $\Si(d)$ be the set of $d$-dimensional cones in $\Si$, and let
$p=|\Si(1)|-3$. Then $X_\Si$ is a GIT quotient
$$
X_\Sigma = \bC^{3+p}\sslash G_\Si =  (\bC^{3+p}- Z_\Si)/G_\Si
$$
where $G_\Si$ is a $p$-dimensional subgroup of $(\bC^*)^{3+p}$ and
$Z_\Si$ is a Zariski closed subset of $\bC^{3+p}$ determined by
the fan $\Si$. If $X_\Si$ is a smooth toric Calabi-Yau 3-fold
then $G_\Si\cong (\bC^*)^p$ and $G_\Si$ acts freely on
$\bC^{3+p}-Z_\Si$. In general we have $(G_\Si)_0\cong (\bC^*)^p$,
where $(G_\Si)_0$ is the connected component of the identity, and
the stabilizers of the $G_\Si$-action on $\bC^{3+p}-Z_\Si$ are at most finite and
generically trivial. The stacky quotient
$$
\cX = [(\bC^{3+p}- Z_\Si)/G_\Si]
$$
is a toric Calabi-Yau 3-orbifold; it is a toric Deligne-Mumford
stack in the sense of Borisov-Chen-Smith \cite{BCS05}.

\subsection{Toric crepant resolution and extended K\"{a}hler classes}
Given a semi-projective simplicial toric Calabi-Yau 3-fold $X_\Si$ which is not
smooth, there exists a subdivision $\Si'$ of $\Si$, such that
$$
X_{\Si'} = (\bC^{3+p+s}-Z_{\Si'})/G_{\Sigma'} \longrightarrow
X_\Si = \big((\bC^{3+p}-Z_\Si)\times (\bC^*)^s\big)/G_{\Si'}.
$$
is a crepant toric resolution, where $X_{\Si'}$ is a smooth
toric Calabi-Yau 3-fold, $s=|\Si'(1)|-|\Si(1)|$, and $G_{\Si'}\cong (\bC^*)^{p+s}$.
$X_{\Si'}$ and $X_\Si$ are GIT quotients of the same $G_{\Si'}$-action
on $\bC^{3+p+s}$ with respect to different stability conditions.

Let $K_{\Si'}\cong U(1)^{p+s}$ be the maximal compact subgroup of $G_{\Si'}\cong (\bC^*)^{p+s}$.
The $G_{\Si'}$-action on $\bC^{3+p+s}$ restricts to a Hamiltonian $K_{\Si'}$-action on
the K\"{a}hler manifold $(\bC^{3+p+s}, \omega_0=\sqrt{-1} \sum_{i=1}^{3+p+s} dz_i \wedge d\bar{z}_i)$,
with moment map $\tmu: \bC^{3+p+s}\to \bR^{p+s}$. There exist two (open) cones
$C$ and $C'$ in $\bR^{p+s}$ such that
\begin{align*}
&\quad\ \tilde{\mu}^{-1}(\vec{r})/K_{\Si'}\\[1ex]
&= \begin{cases}
(\bC^{3+p+s}-Z_{\Si'})/G_{\Sigma'} = X_{\Si'}, & \vec{r}\in C',\\
\big((\bC^{3+p}-Z_{\Si})\times (\bC^*)^s\big)/G_{\Sigma'} = (\bC^{3+p}-Z_{\Si})/G_\Si = X_{\Si}, & \vec{r}\in C
\end{cases}
\end{align*}
$C'\subset \bR^{p+s}= H^2(X_{\Si'};\bR)$ is the K\"{a}hler cone of $X_{\Si'}$ and
$C\subset \bR^{p+s}$ is the extended K\"{a}hler cone of $X_\Si$.

The parameter $\vec{r}\in C$ determines a K\"{a}hler form $\omega(\vec{r})$ on
the toric Calabi-Yau 3-orbifold
$\cX=[(\bC^{3+p}-Z_\Si)/G_\Si]$. The $p+s$ parameters
$\vec{r}=(r_1,\ldots,\allowbreak r_{p+s})$ are extended K\"{a}hler parameters
of $\cX$, where $r_1,\ldots, r_p$ are K\"{a}hler parameters of $\cX$.
The A-model closed string flat coordinates are complexified
extended K\"{a}hler parameters
$$
\tau_a=-r_a+\sqrt{-1}\theta_a,\quad a=1,\ldots, p+s.
$$

\subsection{Toric graphs}
The action of the Calabi-Yau torus $\bT'$ on $\cX$ restricts
to a Hamiltonian $\bT'_\bR$-action on the K\"{a}hler orbifold
$(\cX,\omega(\vec{r}))$, with moment map $\mu':\cX\to M'_\bR=\bR^2$.
The 1-skeleton $\cX^1$ of the toric Calabi-Yau 3-fold $\cX$ is the union of 0-dimensional
and 1-dimensional orbits of the $\bT$-action on $\cX$.
$\Gamma:= \mu'(\cX^1) \subset \bR^2$ is a planar trivalent graph, which
is known as the toric graph of the symplectic toric Calabi-Yau 3-orbifold
$(\cX,\omega(\vec{r}))$. The toric diagram
depends only on the complex structure on $\cX$, where as the toric graph depends
also on the symplectic structure of $\cX$.

\subsection{Aganagic-Vafa Lagrangian branes}
An Aganagic-Vafa Lagrangian brane in a toric Calabi-Yau 3-orbifold $\cX$ is a Lagrangian
sub-orbifold of the form
$$
\cL=[\tL/K_{\Si'}]\subset \cX=[\tmu^{-1}(\vec{r})/K_{\Si'}]
$$
where
\begin{align*}
 \tL = \Bigg\{ &(z_1,\ldots, z_{3+p+s})\in \tmu^{-1}(\vec{r}): \\
&\sum_{i=1}^{3+p+s}\hat{l}_i^1|z_i|^2=c_1, \sum_{i=1}^{3+p+s}\hat{l}_i^2|z_i|^2=c_2,\
 \arg(z_1\cdots z_{3+p+s}) = c_3 \Bigg\},
\end{align*}
$c_1,c_2,c_3$ are constants, and
$$
\sum_{i=1}^{3+p+s} \hat{l}^\alpha_i=0,\quad \alpha=1,2.
$$
The compact 2-torus $\bT'_\bR\cong U(1)^2$ acts
on $\cL$, and $\mu'(\cL)$ is a point on the toric graph $\Gamma=\mu'(\cX^1)$ which
is not a vertex. $\cL$ intersects a unique 1-dimensional $\bT$ orbit $\fl\subset \cX$.
We have $\fl\cong \bC^*\times \cB \bZ_m$ for some positive integer $m$.
When $m=1$, $\cL\cong S^1\times \bC$ is smooth; when $m>1$,
$\cL$ is smooth away from $\cL\cap \fl \cong S^1\times \cB\bZ_m$.

\subsection{Chen-Ruan orbifold cohomology}\label{sec:CR}
Let $U=\bC^{3+p}-Z_\Si$, so that $\cX=[U/G_\Si]$. Given $v\in G_\Si$, let
$U^v = \{ z\in U: v\cdot z= z\}$.
The inertia stack of $\cX$ is
$$
\IX =\bigcup_{v\in \BSi} \cX_v
$$
where $\BSi=\{v\in G_\Si: U^v\neq \emptyset\}$ and $\cX_v =[U^v/G_\Si]$.

We consider cohomology with $\bC$-coefficient. As a graded $\bC$-vector space,
the Chen-Ruan orbifold cohomology \cite{CR04} of $\cX$ is
$$
H^*_\CR (\cX;\bC) = \bigoplus_{v\in \BSi} H^*(\cX_v;\bC)[2\age(v)], \quad \age(v)\in \{0,1,2\}.
$$

Let $\fg :=|\Int(P)\cap N'|$
be the number of lattice points in $\Int(P)$, the interior of the polytope $P$, and
let $\fn:=|\partial P\cap N'|$ be the number of lattice points on $\partial P$, the boundary of the polytope $P$.
Then
\begin{align*}
p &= |\Si(1)|-3= \dim_\bC  H^2(X_\Si;\bC),\\
p+s&= |\Si'(1)|-3= |P \cap N'|-3=  \dim_\bC H^2(X_{\Si'};\bC)=\dim_{\bC} H^2_\CR(\cX;\bC)\\
   &= \fg+\fn-3,\\
\fg &= |\Int(P)\cap N'| =\dim_{\bC} H^4(X_{\Si'}) =\dim_\bC H^4_\CR(\cX;\bC),\\
\chi &= |\Si'(3)|= 2\Area(P)=\dim_\bC H^*(X_{\Si'};\bC)= \dim_\bC H^*_\CR(\cX;\bC)\\
 &= 1+ p+s+\fg = 2\fg-2+\fn.
\end{align*}

\subsection{The mirror curve}

Following the notation in Section \ref{sec:CR}, the convex polytope
$P\subset N'_\bR\cong \bR^2$ defines a polarized toric surface
$(S,L)$, where $S$ is a toric variety and $L$ is an ample line
bundle. We have
$$
\chi (S,L)=h^0(S,L)=|P\cap N'|=3+p+s.
$$
The mirror curve $H(X,Y)$ is given by
$$
H(X,Y)=\sum_{(m,n)\in P\cap N'} a_{m,n} X^m Y^n,\ a_{m,n} \in \bC^*.
$$
So $H(X,Y)\in H^0((\bC^*)^2, \cO_{(\bC^*)^2})$ is the restriction of
a section $s\in H^0(S,L)$. The compactified mirror curve
is $s^{-1}(0)\subset S.$

The element $(t_1,t_2,t_3)\in (\bC^*)^3$ acts on the section $H(X,Y)$
by
$$
H(X,Y) \mapsto t_3H(t_1X, t_2Y).
$$
Modulo this action, the mirror curve is parametrized by $p+s$ elements
$q=(q_1,\dots, q_{p+s})\in (\bC^*)^{p+s}$.  For generic $q$, the mirror curve
$C_q$ is a Riemann surface of genus $\fg$ with $\fn$ punctures, and
the compactified mirror curve $\Cbar_q$ is a smooth hypersurface in
the toric surface $S$. The
Euler characteristic of $C_q$ is
$$
\chi(C_q)=2-2\fg-\fn=-\dim_\bC H^*(\cX;\bC)=-\chi(X).
$$
\subsection{Framings}
The framing $f\in \bZ$ specifies a $1$-dimensional subgroup
$$
\bT_f=\mathrm{ker} (\mathsf{f})\subset \bT',
$$
where the character $\mathsf{f}=\ww'_1- f \ww_2'\in M'=\Hom(N',\bZ).$
It induces a surjective group homomorphism
$$
(\bT')^\vee \cong (\bC^*)^2 \to (\bT_f)^\vee \cong \bC^*, \quad
(X,Y)\mapsto XY^f.
$$
Other than several finite number of choices of $f$, the function
$$
\hat X:=XY^f: C_q\to \bC^*
$$
is holomorphic Morse, i.e. it has simple ramification points. We have
$$
|\Crit(\hat X)|=-\chi(C_q)=2\fg-2+\fn=\dim_\bC H^*_\CR (\cX;\bC).
$$
Around each ramification point $p_0 \in \Crit (\hat X)$, one
writes
$$
\hat x=\hat x(p_0)+\zeta_0^2,
$$
where $\zeta_0$ is the local coordinate around $p_0$. We denote $\check
u^{p_0}=\hat x(p_0)$. It depends on the complex parameter $q$, and is
a \emph{canonical coordinate} of the B-model. For
any $p$ in the neighborhood of $p_0$ we define $\bar p$ by
$$
\zeta_0 (p)=-\zeta_0(\bar p).
$$
We also define a multi-valued holomorphic $1$-form on $C_q$
$$
\Phi=\log Y \frac{d\hat X}{\hat X}.
$$

\section{Gromov-Witten invariants of Toric Calabi-Yau 3-orbifolds}

\subsection{Open Gromov-Witten invariants and A-model open potentials}
Let $\cL$ be an Aganagic-Vafa Lagrangian brane in a toric Calabi-Yau 3-orbifold $\cX$.
Then $\cL$ is homotopic to $S^1\times \cB\bZ_m$, so
$$
H_1(\cL;\bZ)= \pi_1(\cL)= \bZ\times \bZ_m.
$$

Open GW invariants of $(\cX,\cL)$  count holomorphic maps
$$
u: (\Sigma, x_1,\ldots, x_\ell, \partial \Sigma =\coprod_{j=1}^n R_j) \to (\cX,\cL)
$$
where $\Si$ is a bordered Riemann surface with stacky points $x_i=B\bZ_{r_i}$ and
$R_j\cong S^1$ are connected components of $\partial \Sigma$. These invariants depend on
the following data:
\begin{enumerate}
\item the topological type $(g,n)$ of the coarse moduli of the domain, where $g$
is the genus of $\Si$ and $n$ is the number of connected components of $\partial\Sigma$,
\item the degree $\beta'=u_*[\Si]\in H_2(\cX,\cL;\bZ)$,
\item the winding numbers $\mu_1,\ldots,\mu_n \in \bZ$ and the monodromies $k_1,\ldots, k_n \in \bZ_m$,
where $(\mu_j,k_j)= u_*[R_j]\in H_1(\cL;\bZ)=\bZ\times \bZ_m$,
\item the framing $f\in\bZ$  of $\cL$.
\end{enumerate}
We call the pair $(\cL,f)$ a framed Aganagic-Vafa Lagrangian brane.
We write $\vmu=((\mu_1,k_1),\ldots, (\mu_n,k_n))$.  Let $\cM_{g,\ell}(\cX,\cL\mid \beta',\vmu)$ be the moduli space\linebreak
parametrizing maps described above, and let $\Mbar_{g,\ell}(\cX,\cL\mid\beta',\vmu)$
be the partial compactification: we allow the domain $\Si$ to have nodal singularities, and
an orbifold/stacky point on $\Si$ is either a marked point $x_j$ or a node; we require
the map $u$ to be stable in the sense that its automorphism group is finite. Evaluation
at the $i$-th marked point $x_i$ gives a map $\ev_i: \Mbar_{(g,n),\ell}(\cX,\cL\mid\beta',\vec{\mu})\to \IX$.

Given $\gamma_1,\ldots, \gamma_\ell\in H^*_{\CR,\bT'}(\cX;\bC)$, we define
\begin{align*}
\langle \gamma_1,\ldots, \gamma_\ell\rangle_{g,\beta,\vec{\mu}}^{\cX,(\cL,f)}
&:= \int_{[\Mbar_{(g,n),\ell}(\cX,\cL\mid\beta',\vec{\mu})^{\bT_\bR'}]^\vir }
\left.\frac{\prod_{i=1}^\ell \ev_i^*\gamma_i }{e_{\bT'_\bR}(N^\vir)}\right|_{(\bT_f)_\bR}\\
&\in  \bC v^{\sum_{i=1}^\ell \frac{\deg\gamma_i}{2}-1}
\end{align*}
where $v$ is the generator of $H^2(\cB(\bT_f)_\bR ;\bZ)=H^2(\cB U(1);\bZ)\cong \bZ$.

For $\btau=\sum_{a=1}^{p+s} \tau_a e_a \in H^2_\CR(\cX;\bC)$, we define
generating functions $F_{g,n}^{\cX,(\cL,f)}$ of open Gromov-Witten invariants as follows.
\begin{align*}
&\quad\ F_{g,n}^{\cX,(\cL,f)}(Z_1,\ldots, Z_n,\btau)\\
&= \sum_{\beta',\ell\geq 0}\sum_{(\mu_j,k_j)\in \bZ\times \bZ_m}
\frac{\langle \btau^\ell \rangle^{\cX,(\cL,f)}_{g,\beta,(\mu_1,k_1)\cdots,
 (\mu_n, k_n)}}{\ell!} \\
&\quad\   \cdot \otimes_{j=1}^n \Big(Z_j^{\mu_j} (-(-1)^\frac{-k_j}{m})\one'_\frac{-k_j}{m}\Big)
\in H^*_\CR(\cB\bZ_m;\bC)^{\otimes n}
\end{align*}
where $H^*_\CR(\cB\bZ_m;\bC)=\oplus_{k=0}^{m-1} \bC \one'_{\frac{k}{m}}$.

\subsection{Primary closed Gromov-Witten invariants and\\ A-model free energies}

We define genus $g$, degree $\beta$ primary closed Gromov-Witten invariants:
$$
\langle \btau^\ell \rangle^{\cX}_{g,\beta}=
\int_{[\Mbar_{g,\ell}(\cX,\beta)^{\bT_\bR'}]^\vir }
\left.\frac{\prod_{i=1}^\ell \ev_i^*\btau }{e_{\bT'_\bR}(N^\vir)}\right|_{(\bT_f)_\bR}\\
\in  \bC.
$$
This closed Gromov-Witten invariant can be viewed as the case when $n=0$ i.e. there is no boundary on the domain curve. The A-model genus $g$ free energy $F_g^{\cX}$ is
a generating function of primary genus $g$ closed Gromov-Witten invariants.
\begin{eqnarray*}
F_{g}^{\cX}(\btau)
&=& \sum_{\beta,\ell\geq 0}
\frac{\langle \btau^\ell \rangle^{\cX}_{g,\beta}}{\ell!}. \\
\end{eqnarray*}

The BKMP remodeling conjecture builds the mirror symmetry for the open Gromov-Witten potentials $F_{g,n}^{\cX,(\cL,f)}(Z_1,\ldots, Z_n,\btau)$ as well as free energies $F_{g}^{\cX}(\btau)$.

\subsection{Descendant closed Gromov-Witten invariants} Given
$\gamma_1,\ldots,\gamma_n$, we define a generating function of genus $g$, $n$-point
descendant closed Gromov-Witten invariants:
$$
\left\llangle\frac{\gamma_1}{z_1-\psi_1},\ldots,\frac{\gamma_n}{z_n-\psi_n}
\right\rrangle_{g,n}^{\cX}=\sum_{\beta,\ell\geq 0}\frac{1}{\ell!}
\left\langle\frac{\gamma_1}{z_1-\psi_1},\ldots,\frac{\gamma_n}{z_n-\psi_n},
\btau^\ell\right\rangle_{g,\beta}^{\cX},
$$
where $\psi_i=c_1(\bL_i)$ and $\bL_i\to \Mbar_{g,n+\ell}(\cX,\beta)$ is line bundle
whose fiber at moduli point $[u:(C,x_1,\ldots,x_{n+\ell})\to \cX]$ is
the cotangent line $T_{x_i}^*C$ at the $i$-th marked point to (the coarse moduli
space of) the domain curve.

We will state an extension of the remodeling conjecture to higher genus descendent potentials ${\left\llangle\frac{\gamma_1}{z_1-\psi_1},\ldots,\frac{\gamma_n}{z_n-\psi_n}
\right\rrangle_{g,n}^{\cX}}$.

\section{Eynard-Orantin invariants of the mirror curve}

\subsection{Fundamental normalized differential of the second kind}\label{sec:Bergmann}

In this subsection, we recall the definition of the fundamental normalized differential
of the second kind $B(p_1,p_2)$ (see e.g. \cite{Fay}) for a general compact
Riemann surface $\Cbar$.

Let $\Cbar$ be a compact Riemann surface of genus $\fg$. When $\fg>0$, let
$A_1,B_1,\allowbreak \ldots, A_{\fg}, B_{\fg}$ be a symplectic basis of $(H_1(\Cbar;\bC), \cdot)$:
$$
A_i\cdot A_j = B_i\cdot B_j=0,\quad A_i\cdot B_j = -B_j\cdot A_i = \delta_{ij}
$$
where $\cdot$ is the intersection pairing. For our purpose, we need
to consider $H_1(\Cbar;\bC)$ instead of the integral first homology group $H_1(\bar{C};\bZ)$.
We assume that the Lagrangian subspace $\bigoplus_{i=1}^{\fg}\bC A_i$ of $H_1(\Cbar;\bC)$
is transversal to the Lagrangian subspace
$$
H^{1,0}(\Cbar)^{\perp}:= \{ \gamma\in H_1(\Cbar;\bC): \langle \theta,\gamma\rangle=0
\quad \forall \theta \in H^{1,0}(\Cbar)\}
$$
where $\langle \  , \ \rangle: H^1(\Cbar;\bC)\times H_1(\Cbar;\bC)$ is the
natural pairing; this assumption holds when $A_1,\ldots, A_{\fg}\in H_1(\Cbar;\bR)$.

\enlargethispage{1em}
The fundamental normalized differential of the second kind $B(p_1,p_2)$  on $\bar{C}$ is characterized by the following properties:
\begin{enumerate}
\item $B(p_1,p_2)$ is a bilinear symmetric meromorphic differential on $\Cbar_q\times\Cbar_q$.
\item $B(p_1,p_2)$ is holomorphic everywhere except for a double pole along the
diagonal.  If $p_1,p_2$ have local coordinates $z_1,z_2$ in an open neighborhood $U$
of $p\in \Cbar_q$ then
$$
B(p_1,p_2)=\left(\frac{1}{(z_1-z_2)^2} + a(z_1,z_2)\right) dz_1 dz_2
$$
where $a(z_1,z_2)$ is holomorphic on $U\times U$ and symmetric in $z_1,z_2$.
\item $\displaystyle{\int_{p_1\in A_i}B(p_1,p_2)}=0$, $i=1,\ldots,\fg$.
\end{enumerate}
In fact, we do not need a particular choice of A-cycles. We just need to
specify a Lagrangian subspace of  $(H_1(\Cbar;\bC),\cdot)$ transversal
to the Lagrangian subspace $H^{1,0}(\Cbar;\bC)^\perp \subset H_1(\Cbar;\bC)$ such that
the period of $B(p_1,p_2)$ along any element in this subspace is zero.

The fundamental differential $B(p_1,p_2)$ also satisfies the following properties:
\begin{enumerate}
\item[(4)] If $f$ is a meromorphic function on $\Cbar$ then
$$df(p_1)= \displaystyle{ \Res_{p_2\to p_1}B(p_1,p_2) f(p_2)}.$$
\item[(5)] $\displaystyle{\int_{p_1\in B_i} B(p_1,p_2) = 2\pi\sqrt{-1}\omega_i(p_2)}$,
where $\omega_i$ is the unique holomorphic 1-form on $\Cbar$ such that
$\int_{A_j}\omega_i =\delta_{ij}$.
\end{enumerate}

\subsection{Choice of A-cycles on the compactified mirror curve}
The mirror theorem for semi-projective toric orbifolds \cite{CCIT}
relates the $1$-primary $1$-descendant function (the $J$-function)
$$
\left\llangle 1, \frac{\phi_a}{z-\psi}\right\rrangle^\cX_{0,2}\phi^a
$$
to certain hypergeometric $I$-function $I^\cX(q,z)$ under the mirror map
$$
\tau_a =\frac{1}{2\pi\sqrt{-1}}\int_{A_a}\Phi =
\begin{cases}
\log q_a + h_a(q), & a=1,\ldots, p\\
q_a (1+ h_a(q)), & a=p+1,\ldots, p+s.
\end{cases}
$$
which as the prescribed leading term behavior (all $h_a(q)$ are power
series in~$q$).

It is a well-known fact that these mirror maps are
given by such period
integrals where $A_a\in H_1(C_q;\bC)$.
The inclusion $C_q\hookrightarrow \Cbar_q$ induces a surjective
group homomorphism  $H_1(C_q;\bC)\cong \bC^{2\fg+\fn-1} \to H_1(\Cbar_q)\cong \bC^{2\fg}$ where the kernel is generated by the $\fn$ loops around the $\fn$ points in
$\Cbar_q\setminus C_q$; each of these $\fn$ loops is contractible
in $\Cbar_q$, and the sum of these $\fn$ loops is homologous to zero in $C_q$. The images of $A_a\in H_1(C_q;\bC)$ in $H_1(\Cbar_q;\bC)$ span a Lagrangian subspace $L_A\subset H_1(\Cbar_q;\bC)$ transversal to the Lagrangian subspace
$H^{1,0}(\Cbar_q)^\perp \subset H_1(\Cbar_q;\bC)$.
We use the Lagrangian subspace $L_A$ to define our fundamental normalized differential of the second kind
$B(p_1,p_2)$ for the purpose of constructing higher genus B-model
invariants.

\subsection{The Eynard-Orantin topological recursion}

We use the fundamental differential $B$ prescribed above to run the
Eynard-Orantin topological recursion. It starts with two initial data
(unstable cases)
$$
\omega_{0,1}=0,\quad \omega_{0,2}=B.
$$
The stable cases ($2g-2+n>0$) are defined recursively by the Eynard-Orantin
topological recursion:
\begin{align*}
&\omega_{g,n}(p_1,\ldots,p_n)\\
:=& \sum_{p_0\in\Crit(\hX)}
\Res_{p\to p_0}\frac{\int_{\xi=p}^{\bar{p}}B(p_n,\xi)}{2 (\Phi(p)-\Phi(\bar{p}))}
\Biggl(\omega_{g-1,n+1}(p,\bar{p},p_1,\ldots, p_{n-1})\\
& +\sum_{\substack{g_1+g_2=g \\ I\sqcup J=\{1,\ldots,n-1\} }} \omega_{g_1,|I|+1}(p,p_I)
\omega_{g_2,|J|+1}(\bar{p},p_J)\Biggr).
\end{align*}
The resulting $\omega_{g,n}$ for $2g-2+n>0$ is a symmetric meromorphic form on
$(\Cbar_q)^n$. They are holomorphic on $(\Cbar_q\setminus\Crit(\hX))^n$
and satisfy the following properties:
\begin{enumerate}
\item For any $j\in \{1,\ldots,n\}$ and  any $p_0\in \Crit(\hX)$,
$$
\Res_{p_j\to p_0} \omega_{g,n}(p_1,\ldots,p_n) =0.
$$
\item For any $j\in \{1,\ldots,n\}$ and  any $i\in\{1,\ldots,\fg\}$,
$$
\int_{p_j\in A_i} \omega_{g,n}(p_1,\ldots,p_n) =0.
$$
\end{enumerate}

\subsection{B-model open potentials}

For $\ell\in \bZ_m=\bZ/m\bZ$, let
$$
\psi_\ell :=\frac{1}{m}\sum_{k=0}^{m-1} e^\frac{2\pi\sqrt{-1}k\ell}{m} \one'_\frac{k}{m}.
$$
Then $\{ \psi_\ell: \ell=0,1\ldots, m-1\}$ is a canonical basis of  $H^*_\CR(\cB\bZ_m;\bC)$.

Recall that $\cL$ intersects a unique 1-dimensional orbit $\fl$ of the $\bT$-action on $\cX$.
We assume that the closure $\bar{\fl}$ of $\fl$ in $\cX$ is non-compact, so that $\cL$ is
an ``outer'' brane. Then the 2-dimensional cone associated to $\bar{\fl}$ corresponds
an edge $e$ on the boundary of the polytope $P$, and $|e\cap N'|=m+1$. Let
$D\subset S$ be the torus invariant divisor associated to the edge $e$. For generic $q$,
the compactified mirror curve $\Cbar_q$ intersects $D$ transversally at $m$ points
$\bar{p}_0 ,\ldots,\bar{p}_{m-1}$. For $\ell\in \{0,1,\ldots,m-1\}$, there
exist open neighborhoods $U_\ell$ of $\bar{p}_\ell$ in the compactified mirror
curve $\Cbar_q$ and $U$ of $0$ in $\bP^1=\bC^*\cup \{0,\infty\}$ such that
$\hat{X}|_{U_\ell}:U_\ell\to U$ is biholomorphic. Let
$\rho_\ell:= (\hat{X}|_{U_\ell})^{-1}: U\to U_\ell$. We define
B-model topological open string partition functions as follows.
\begin{enumerate}
\item disk invariants
$$
\check{F}_{0,1}(q;X): =\sum_{\ell\in \bZ_m} \int_0^X \left( \big(\log Y(\rho_\ell(X'))-\log Y(\bar{p}_\ell)\big)
\frac{dX'}{X'}\right) \psi_\ell
$$
which take values in $H^*_\CR(\cB\bZ_m;\bC)$.

\item annulus invariants
\begin{align*}
&\quad\ \check{F}_{0,2}(q;X_1,X_2)\\
&:= \sum_{\ell_1,\ell_2\in \bZ_m}
\int_0^{X_1}\int_0^{X_2} \left((\rho_{\ell_1}\times \rho_{\ell_2})^*\omega_{0,2}-\frac{dX'_1 dX'_2}{(X'_1-X'_2)^2} \right) \psi_{\ell_1}
\otimes \psi_{\ell_2}
\end{align*}
which take values in $H^*_\CR(\cB\bZ_m;\bC)^{\otimes 2}$.

\item $2g-2+n>0$
\begin{align*}
&\quad\ \check{F}_{g,n}(q;X_1,\ldots,X_n)\\
&:=\sum_{\ell_1,\cdots,\ell_n\in\bZ_m}\int_0^{X_1}\cdots \int_0^{X_n}
(\rho_{\ell_1}\times \cdots \times \rho_{\ell_n})^*\omega_{g,n}
\psi_{\ell_1}\otimes\cdots \otimes \psi_{\ell_n}
\end{align*}
which take values in $H^*_\CR(\cB\bZ_m;\bC)^{\otimes n}$.
\end{enumerate}
Each of the $m^n$ components of $\check{F}_{g,n}(q;X_1,\ldots,X_n)$  is
a power series in $q_1,\ldots,q_{p+s},X_1,\ldots, X_n$ which converges in an open neighborhood of the origin.

\subsection{B-model free energies}

For $g\ge 2$, the B-model free energy is defined as
$$
\check F_g(q) =\frac{1}{2g-2}\sum_{p_0\in \Crit{\hX}}\Res_{p\to p_0} \omega_{g,1}(p) \widetilde \Phi(p),
$$
where
$$
d \widetilde \Phi= \Phi.
$$
Notice that the function $\widetilde \Phi$ locally defined around each
critical point of $\hX$ has some ambiguities, since $\Phi$ is
multi-valued, and $\widetilde \Phi$ is determined by $\Phi$ up to a
constant. However, the residue is well-defined since is does not
depend on these ambuities.

For $g=1$, the free energy is defined up to a constant
$$
\check F_1(q)=-\frac{1}{2}\log \tau_B-\frac{1}{24}\sum_{p_0\in
  \Crit{\hX}} \log h^{p_0}_1.
$$
Here the Bergmann $\tau$-function $\tau_B$ is defined up to a constant by
$$
 \frac{\partial \log \tau_B}{\partial \check
   u^{p_0}}=\Res_{p\to p_0}\frac{B(p,\bar p)}{d\hat x(p)}.
$$

When $g=0$, the prepotential $F_0$ is characterized by
$$
\frac{\partial \check F_0}{\partial \tau_a}=\int_{p\in B_a} \Phi(p).
$$
Notice that since $\Phi$ is a multi-valued differential form, and it
satisfies the following
$$
\frac{\partial \Phi(p)}{\partial \tau_a}=\int_{p'\in B_a} \omega_{0,2}(p,p').
$$
The prepotential $\check F_0$ defined this way is only determined up to
a quadratic polynomial in $\tau_a$.

\section{All genus open-closed mirror symmetry}

In this section, $(\cL,f)$ is an outer Aganagic-Vafa Lagrangian brane in $\cX$,
so that the closure of $\fl=\bC^*\times \cB \bZ_m$ contains
a unique $\bT$ fixed point. Let $G$ be the stabilizer of this fixed point. Then $G$  is a finite abelian group which
contains $\bZ_m$ as a subgroup. When $\cX$ is smooth, we have $m=1$ and $G$ is trivial.

\subsection{All genus open-closed mirror symmetry:\\ the remodeling conjecture}

\begin{conjecture}[Bouchard-Klemm-Mari\~{n}o-Pasqetti \cite{BKMP09, BKMP10}] \label{remodel}
$$
\check{F}_{g,n}(q;X_1,\ldots,X_n) =(-1)^{g-1+n}|G|^n F_{g,n}^{\cX,(\cL,f)}(\btau; Z_1,\ldots,Z_n)
$$
where
$(q,X_j)$ and $(\btau,Z_j)$ are related by the open-closed mirror map:
\begin{align*}
\tau_a =\frac{1}{2\pi\sqrt{-1}}\int_{A_a}\Phi &=
\begin{cases}
\log q_a + h_a(q), & a=1\ldots, p\\
q_a (1+ h_a(q)), & a=p+1,\ldots,p+s
\end{cases}\\[1ex]
\log Z_j &=  \log X_j + h_0(q)
\end{align*}
where $h_0(q), h_1(q),\ldots, h_{p+s}(q)$ are explicit power series in $q$ convergent in
a neighborhood of the origin in $\bC^{p+s}$. Notice that when $n=0$, this is a
statement about closed Gromov-Witten mirror symmetry (and the
right-hand side does not depend on $(\cL,f)$). When $(g,n)=(1,0)$ and
$(0,0)$, the statement is valid up to a constant and a quadratic
polynomial in $\tau_a$, respectively.
\end{conjecture}
Indeed, the above statement is more general than the original conjecture
in \cite{BKMP09, BKMP10}, where they conjecture about non-gerby branes
(the $m=1$ case).

Conjecture \ref{remodel} was proved when $\cX=\bC^3$ independently
by L. Chen \cite{Ch09} and J. Zhou \cite{Zh09}.
In 2012, Eynard-Orantin provided a proof of the BKMP remodeling conjecture for all symplectic smooth toric Calabi-Yau 3-folds \cite{EO12}.
In the orbifold case, the authors prove Conjecture \ref{remodel} first for
affine toric Calabi-Yau 3-orbifolds \cite{FLZ13} and later for
all semi-projective toric Calabi-Yau 3-orbifolds \cite{FLZ15}.

We now give a brief outline of the proof of Conjecture \ref{remodel} in \cite{FLZ15}.
Givental proved a quantization formula for total descendant
potential of equivariant GW theory of GKM manifolds \cite{Gi01, Gi, Gi04}. (See also the book by
Lee-Pandharipande \cite{LP}.) The third author generalized this formula
to GKM orbifolds \cite{Zo14}.  The quantization formula is equivalent to a graph sum formula of the total descendant
potential, which implies a graph sum formula
$$
F_{g,n}^{\cX,(\cL,f)}= \sum_{\vGa\in G_{g,n}}\frac{w_A(\vGa)}{|\Aut(\vGa)|},
$$
where $G_{g,n}$ is a certain set of decorated stable graphs.
The unique solution $\{\omega_{g,n}\}$ to the Eynard-Orantin topological recursion can be expressed as a sum over graphs \cite{KO, E11, E14, DOSS}.
We expand the graph sum formula in \cite[Theorem 3.7]{DOSS}
(which is equivalent to \cite[Theorem 5.1]{E11})
 at punctures $\{ \bar{p}_\ell:\ell\in \bZ_m\}$, and obtain
a graph sum formula
$$
\check{F}_{g,n} = \sum_{\vGa\in G_{g,n}}\frac{w_B(\vGa)}{|\Aut(\vGa)|}.
$$
Finally, we use the genus-zero mirror theorem for smooth toric DM stacks \cite{CCIT}
to prove
$$
w_B(\vGa)= (-1)^{g-1+n} |G|^n w_A(\vGa)
$$
for all decorated graphs $\vGa$.

\subsection{Descendant version of the all genus mirror symmetry}

Iritani \cite{Ir09} studies the oscillatory integral and shows the
following
$$
\int_{\SYZ(\cF)} e^{\frac{\widetilde W}{z}}\Omega=\left\llangle \frac{\kappa(\cF)}{z-\psi}\right\rrangle_{0,1},
$$
where $\cF$ is a
$\bT_f$-equivariant coherent sheaf on $\cX$. Here the $\SYZ$ is the SYZ T-dual functor,
which takes a $\bT_f$-equivariant coherent sheaf on $\cX$ and produces a Lagrangian
brane in $(\bC^*)^3$.\footnote{Iritani \cite{Ir09} does not
  explicitly states this identity under the SYZ transform, but instead
  he matches the cases $\cF=\cO_{\cX}$ and a skyscraper sheaf. He then
applies the monodromy to $\cO_{\cX}$ to obtain other line bundles
on $\cX$. These sheaves generate the K-theory group.} The \emph{equivariantly perturbed} superpotential $W$
is given by
$$
\tW=W-\log X - f\log Y.
$$

Let $(\bT_f)_\bR\cong U(1)$ be the maximal torus of $\bT_f\cong\bC^*$, and
let $\mu_{(\bT_f)_\bR}:\cX\to \bR$ be the moment map of the Hamiltonian
$(\bT_f)_\bR$-action on $(\cX,\omega)$.
We say a $\bT_f$-equivariant coherent sheaf $\cF$ on $\cX$ is admissible if
(i) $\mu_{(\bT_f)_\bR}(\mathrm{supp}(\cF))\subset \bR$ is bounded below,
and (ii) the Lagrangian brane $\SYZ(\cF)$ reduces to a cycle
$\gamma(\cF)$ on the mirror curve $C_q$, while the oscillatory integral
could be done on the curve
$$
\int_{\mathrm{SYZ}(L)}e^{\frac{\widetilde W}{z}} \frac{dXdYdZ}{XYZ}=
\int_{\gamma(\cF)} e^{\frac{\hat x}{z}} ydx.
$$
Condition (i) implies that $\hat{x}$ is bounded below on
$\gamma(\cF)$, so the integral on the RHS converges when $z\in (-\infty,0)$.

Using this result and combining with the remodeling conjecture, we
have
\begin{theorem}[Descendant version of the all genus mirror symmetry
  for $\cX$]
$$
\int_{\gamma(L_1)\times \dots\times \gamma(L_n)} e^{\frac{\hat
    x_1}{z_1}+\dots +\frac{\hat x_n}{z_n}}\omega_{g,n}=\left\llangle
\frac{\kappa(L_1)}{z_1-\psi_1} \dots \frac{\kappa(L_n)}{z_n-\psi_n}
\right\rrangle_{g,n}.
$$
\end{theorem}
To obtain this theorem, one observes that when integrating
$\omega_{g,n}$ we are simply integrating the leaf terms of
$w_A(\vGa)$, since only leaf terms are forms while all other graph
component contributions are scalars. The genus $0$ oscillatory
integral theorem from \cite{Ir09} turns these leafs into genus $0$ descendants, and the
graph becomes precisely the graph for higher genus descendant potentials.

\end{document}